\documentclass[english,reqno, 11pt]{amsart}
\usepackage[T1]{fontenc}
\usepackage[latin9]{inputenc}
\usepackage{amsthm}
\usepackage{amssymb}
\usepackage{amsmath}
\usepackage{mathabx}
\usepackage{xcolor}
\usepackage{verbatim, color}
\usepackage{bbm}
\usepackage{float}
\usepackage{dsfont, graphicx}
\usepackage{epstopdf}
\usepackage[hidelinks]{hyperref}
\usepackage[T1]{fontenc}
\usepackage{amssymb,amsmath, amsthm, amsfonts}
\usepackage{graphicx}

\usepackage[hidelinks]{hyperref}

\usepackage{thmtools,thm-restate}

\makeatletter
\numberwithin{equation}{section}
\numberwithin{figure}{section}
\theoremstyle{plain}
\newtheorem{thm}{Theorem}[section]

\newtheorem{lem}[thm]{Lemma}

\newtheorem{rem}[thm]{Remark}

  \newcounter{casectr}
  
\theoremstyle{definition}

\theoremstyle{remark}

\newcommand{\RRR}{\mathbb{R}}

\newcommand{\qb}{\tilde{Q}_{b}}

\def\<{\langle}
\def\>{\rangle}

\title{A note on log-log blow up solutions for stochastic nonlinear Schr\"odinger equations}

\begin{document}
\begin{center}
\author{Chenjie Fan}
\address[C. Fan]{Academy of Mathematics and Systems Science, Chinese Academy of Sciences, Beijing , China}
\author{Yiming Su}
\address[Y. Su]{School of science,
Zhejiang University of Technology, Hangzhou, China.}
\author{Deng Zhang}
\address[D. Zhang]{School of mathematical sciences,
Shanghai Jiao Tong University,  Shanghai, China.}
\end{center}
\maketitle
\begin{abstract}
In this short note, we present a construction for the log-log blow up solutions to focusing
mass-critical stochastic nonlinear Schr\"oidnger equations  with multiplicative noises.
The solution is understood in the sense of controlled rough path as in \cite{SZ20}.
\end{abstract}

\section{Introduction}
In this article, we investigate the  stochastic mass critical focusing nonlinear Schr\"odinger equation 
\begin{equation}\label{eq: rough}
\begin{cases} 
      dX = i\Delta Xdt + i|X|^{\frac 4d}X dt - \mu X dt + iXdW(t),  \\
      X(0)= X_0 \in H^1(\RRR^d). 
\end{cases}
\end{equation}
Here, $W$ is the a Wiener process colored in space and white in time.
$$W(t,x)=\sum_{k=1}^{N} \phi_k(x)B_k(t),\ \ x\in \RRR^d,\ t\geq 0,$$
where
$\phi_k$ is real valued, and
$B_k$ are the standard $N$-dimensional real valued Brownian motions
on a stochastic basis $(\Omega, \mathcal{F}, \{\mathcal{F}_t\}, \mathbb{P})$,
$1\leq k\leq N <\infty$. And $\mu= \frac{1}{2} \sum_{k=1}^N  \phi_k ^2$. 
The last term in \eqref{eq: rough} is taken in a certain sense of
controlled rough path.

Throughout this paper,
we only consider very nice noise for simplicity and concreteness.
More precisely,
we assume
\begin{equation}
  \phi_{k} \text{ are Schwartz functions},\ \ \forall 1\leq k\leq N.
\end{equation}

Since $W$ is real valued, the noise is of conservative type. One physical significance is that \eqref{eq: rough} admits a pathwise mass conservation law.

Another physical application can be found in the study of monolayer Scheibe aggregates.
In crystals  the noise corresponds to the scattering of excitons
by phonons due to the thermal vibrations of molecules,
and the noise effects on the dynamics of the two dimensional
nonlinear Schr\"odinger equation was studied
(see \cite{BCIR94}, \cite{BCIRG95}).
See also \cite{RGBC95} for the noise effect on collapse in the one dimensional case.
For more general noise,
including the non-conservative case,
we refer to \cite{BG09} for the application in the quantumn measurement,
where the noise represents the output of continuous measurement.

The well posedness of stochastic nonlinear Schr\"odinger equations
is extensively studied in literature.
See, e.g.,
\cite{BD99,BD03,BM13,BRZ14,BRZ16,H18,FX18.1,FX18.2,Z18}.

We are interested in the blow up dynamic of \eqref{eq: rough}. 

Parallel to the deterministic  mass critical NLS, \cite{dodson2015global}, \cite{weinstein1983nonlinear}, 
\begin{equation}\label{eq: nls}
iu_{t}+\Delta u =-|u|^{4/d}u,
\end{equation}
it is known that equation \eqref{eq: rough}
admits a global flow if its corresponding mass
is below that of the ground state,
i.e., $\|X_0\|_{L^2} < \|Q\|_{L^2}$.
Here, the ground state $Q$ is the positive radial solution to the elliptic equation
\begin{align} \label{eq: Q}
     \Delta Q - Q + |Q|^{\frac 4d} Q  =0.
\end{align}

However, when mass is large, 
the solutions may formalize singularity.
We refer to \cite{BD05} for the characterization of a blow-up region
corresponding to sufficiently negative energy,
based on the stochastic version of virial evolution, \cite{glassey1977blowing}.
Quite interestingly,
in the focusing mass-supsercritical case,
the conservative noise can
accelerate blow-up with positive probability (see \cite{BD05}).
Moreover, in the recent work \cite{SZ20}
minimal mass blow-up solutions are constructed in the stochastic setting,
which actually exhibit pseudo-conformal blow up rate near the blow up time.

It should be also mentioned that,
several numerical results have been studied on the noise effects on blow-up dynamics in the stochastic case.
See, e.g., \cite{BDM01,DM02,DM02.2}.
In particular,
we refer to the recent works \cite{MRRY20,MRRY20.0}
for the study of noise effects on the log-log blow-up dynamics.

The Log-log blow up solution is one of the most well understood blow up solutions to \eqref{eq: nls}. As aforementioned, no solution with mass below the mass of ground state can blow up. And, it is shown by \cite{merle1993determination} that all $H^{1}$ finite time blow up solutions with mass $\|Q\|_{2}^{2}$ are essentially of form $S(t,x):=\frac{1}{|t|^{d/2}}Q(\frac{x}{t})e^{-\frac{i}{t}+\frac{i|x|^{2}}{4t}}$. Next, one may consider initial data with mass just above the ground state,
\begin{equation}
\|Q\|_{2}<\|u_{0}\|_{2}<\|Q\|_{2}+\alpha, \text{ where  } \alpha \text{ small}.
\end{equation}
It has been shown in a series of work by Merle and Rapha\"el, \cite{merle2005blow},\cite{merle2003sharp},\cite{merle2004universality},\cite{merle2006sharp} that all negative energy solutions to \eqref{eq: nls} will blow up in finite time $T$ according to log-log law. Log-log blow up solutions were numerically  observed in \cite{landman1988rate}, and first mathematically constructed in \cite{perelman2001blow}. Note that log-log blow up dynamic is stable under $H^{1}$ perturbation, \cite{raphael2005stability}. Indeed it is stable under $H^{s}, s>0$ perturbation \cite{colliander2009rough}, and it is stable under certain randomized $L^{2}$ perturbation, \cite{fan2020construction}. It is also highly localized \cite{planchon2007existence}, \cite{raphael2009standing}, \cite{holmer2012blow}.

The interest of this note is to take advantage of the nice property of log-log blow up solutions to do construction for the stochastic model.

There are two main simple but useful observations.

\begin{itemize}
\item The rough path formulation of the solution to \eqref{eq: rough} in \cite{SZ20} seems more flexible compared to the usual Ito formulation, and in particular, one can apply a rescaling transformation and reduced the stochastic PDE \eqref{eq: rough} into a deterministic  model with random coefficients, \eqref{eq: RNLS}.  It enables one to focus on a certain collection of paths rather than all the paths.
\item The construction of log-log blow up solutions to \eqref{eq: RNLS} can be seen as an easier task compared to \cite{raphael2006existence},\cite{raphael2009standing}, \cite{colliander2009rough}, thus one would well expect such a construction is doable. Indeed, the only extra ingredient we need is to prove an energy estimate, Lemma \ref{lem: enegyestimate}. With such an estimate, the construction will follow from the robust bootstrap scheme in \cite{planchon2007existence}.
\end{itemize}

See more details in Section \ref{sec: seccon}. 

Our main result is 
\begin{thm}  \label{Thm-X-loglog}
Consider equation \eqref{eq: rough} with $d=1,2$.
Then, there exists an initial datum $X_0 \in H^1$
such that,
with high probability,
the corresponding solution $X$ to \eqref{eq: rough}
blows up in finite time $T=T_{\omega}$ according to  the log-log law in the sense that,
there exist parameters $(x(t), \gamma(t), \lambda(t)) \in C^1((0,T); \RRR^d \times \RRR \times \RRR^{+})$
such that
\begin{equation}
   e^{-iW(t,x)} X(t,x)=\frac{1}{\lambda^{d/2}(t)}(Q+\epsilon)(t, \frac{x-x(t)}{\lambda(t)})e^{i\gamma(t)} ,\ \ t\in (0,T),
\end{equation}
with
\begin{equation}
\lambda(t)^{-1}\sim \sqrt{\frac{\ln |\ln (T-t)|}{T-t}}, \text{ and } \int |\nabla \epsilon|^{2}+|\epsilon|^{2}e^{-|y|} dy \xrightarrow {t\rightarrow T} 0.
\end{equation}
\end{thm}
\begin{rem}
For those familiar with the property of Gaussian and the blow up solution construction, if one can prove such result with positive probability, one can prove such a result for any probability close to 1.  \textbf{However, Theorem \ref{Thm-X-loglog} is far from proving the existence of an (deterministic) initial data, which will blow up according to log-log law with probability 1.} Though, it is not hard to use our construction to find random initial data in $L_{\omega}^{\infty}L_{x}^{2}$, such that it blows up according to the  log-log law almost surely.
\end{rem}
\begin{rem}
For those only interested in the blow up rate, since our $W$ is regular in the space variable (it is a colored noise), when $t$ approaches the blow up time, one has 
\begin{equation}
\|u\|_{H^{1}}\sim \|X\|_{H^{1}}\sim \sqrt{\frac{\ln |\ln (T-t)|}{T-t}}, \text{ as } t\rightarrow T.
\end{equation}
\end{rem}
\begin{rem}
Observe as $t$ approaches blow up time $T^{\omega}$, the leading profile of $X$ is $Qe^{iW(t, \lambda(t)x+x(t))}$. It is known $x(t)$ has a a limit $x^{*}$ as $t\rightarrow T^{\omega}$ in those kinds of blow up. Thus, asymptotically, $Qe^{iW(t, \lambda(t)x+x(t)}$ will converge to $Qe^{iW(T_{\omega},x^{*})}$ in $H^{1}$. This is consistent  with  Conjecture 1 in \cite{MRRY20}.
\end{rem}
\subsection{Notation}
We say $A\lesssim B$ if there is a constant $C$, so that $A\leq CB$. The constant may change line by line. If $A\lesssim B$ and $B\lesssim A$, we say $A\sim B$.

\subsection{Acknowledgment}

C.Fan was supported by a start up funding from AMSS. Y. Su was supported in NSFC (No. 11601482). D. Zhang was supported by NSFC (No. 11871337).

\section{Construction of the log-log blow up solutions}\label{sec: seccon}
\subsection{Step 0: What do we mean by a solution?}
As in \cite{SZ20}, we use a re-scaling transformation
$u=e^{-iW}X$ to reduce the original stochastic equation \eqref{eq: rough}
to a nonlinear Schrodinger equation with random coefficients,
\begin{equation}  \label{eq: RNLS}
\begin{cases}
 i\partial_t u+ e^{-iW}\Delta(e^{iW}u)+|u|^{\frac{4}{d}}u =0,   \\
u(0,x)=X_{0}   \nonumber
\end{cases}
\end{equation}
Note that we have $e^{-iW}\Delta(e^{iW}u) = \Delta u+  \sum_{j}i\beta_{j}\partial_{j} u+ c u$ with
 
\begin{equation}
\begin{aligned}
  i\beta_{j}(t,x)&=2 \partial_{j} iW(t,x) = 2i \sum_{k=1}^N \partial_{j}\phi_k(x) B_k(t),   \\
  c (t,x)  &= - \sum_{j=1}^d (\sum_{k=1}^N \partial_j \phi_k(x) B_k(t))^2
         + i\sum_{k=1}^N \Delta \phi_k(x) B_k(t).    \end{aligned}   
\end{equation}
Such transformation can be viewed as a Doss-Sussman type transformation
in the infinite dimensional space.

The solutions to equation \eqref{eq: RNLS} is understood in the usual mild sense.
It is known (see \cite{SZ20}) that
equation \eqref{eq: RNLS} is local well-posed
and generates continuous flows in the space $H^{1}$,
and the corresponding solutions blow up if and only if their $H^{1}$-norms blow up.

Furthermore, as illustrated in \cite[Theorem 2.13]{SZ20},
if $u$ solves  \eqref{eq: RNLS} on some random\footnote{Since the coefficients of \eqref{eq: RNLS} are random, thus the time interval where $u$ is wellposed is also random} time interval $[0,\tau^*)$,
then $X:= e^W u$ solves \eqref{eq: rough} on $[0,\tau^*)$
in certain sense of controlled rough path sense.   Briefly speaking, this means for any $\phi\in C_c^\infty$,
$t \mapsto \<X(t), \phi\>$ is continuous on $[0,\tau^*)$
and for any $0<s<t<\tau^*$,
\begin{equation}
\begin{aligned}
   \<X(t)-X(s), \phi\>
   - \int_s^t \<i X, \Delta\phi\>  + \<i|X|^{\frac 4d} X, \phi\>  - \< \mu X, \phi\> dr
   = \sum_{k=1}^N \int_s^t \<i\phi_k X, \phi\> dB_k(r).
\end{aligned}
\end{equation}
Here the integral $\int_s^t \<i\phi_k X, \varphi\> d B_k(r)$
is taken as a (pathwise) rough integration rather than a classical stochastic Ito's integration.
We refer to Definition 2.1 in \cite{SZ20} for more details, and we will not explicitly use it in the current article. See also  Theorem 2.3 in \cite{SZ20}.

We refer to seminal work \cite{lyons1998differential}, \cite{gubinelli2004controlling} for the notion of rough path and controlled rough path, see also recent textbook \cite{FH14}.

This pathwise formulation of solutions seems more flexible in the blow up analysis compared to the It\^o formulation,
since it enables one to focus on a certain collection of paths rather than all the paths.

Furthermore, it enables one to  work on \eqref{eq: RNLS}, which is with random coefficients but of deterministic nature, and which can be modelled by a variation of classical NLS, which was extensively studied in the literature.

In the rest of the article, we focus on the construction of log-log blow up dynamic
of equation \eqref{eq: RNLS}, as stated in Theorem \ref{thm: tech} below.  

\begin{thm}\label{thm: tech}
Let $d=1,2$ and $u$ solve \eqref{eq: workingnls} with the initial data $u_{0}$ given by \eqref{eq: initialdata}.
Then,
with high probability,
$u$ blows up in finite time $T$ according to log-log law in the sense that,
there exist parameters
$(x(t), \gamma(t), \lambda(t)) \in C^1((0,T); \RRR^d \times \RRR \times \RRR^{+})$,
such that
\begin{equation}
u(t,x)=\frac{1}{\lambda^{d/2}(t)}(Q+\epsilon)(t, \frac{x-x(t)}{\lambda(t)})e^{i\gamma(t)} ,\ \ t\in (0,T),
\end{equation}
with
\begin{equation}
\lambda(t)^{-1}\sim \sqrt{\frac{\ln |\ln T-t|}{T-t}}, \text{ and } \int |\nabla \epsilon|^{2}+|\epsilon|^{2}e^{-|y|} dy \xrightarrow {t\rightarrow T} 0.
\end{equation}
\end{thm}
\begin{rem}\label{rem: hp}
The notion of high probability may be understood in the following way. The initial data we described in \eqref{eq: initialdata} has a parameter $\lambda_{0}$, and the probability that such solution blows up according to log-log law is at least $1-e^{-\lambda_{0}^{-\delta}}$, $\delta>0$, this is standard in Gaussian type analysis. See also remark \ref{rem: furtherremark}.
\end{rem}

\subsection{Step 1: Starting point}

Let us go back to \eqref{eq: RNLS} and rewrite it as
\begin{equation}\label{eq: workingnls}
iu_{t}+\Delta u+\sum_{j}i\beta_{j}\partial_{j}u+cu + |u|^{4/d}u =0.
\end{equation}
Up to dropping a set of small probability,
we focus on the paths satisfying
\begin{equation}\label{eq: assumptionsforbandc}
\sup_{t\in [0,1]}\sum_{j}\|\beta_{j}\|_{M}+\|c\|_{M}\lesssim 1.
\end{equation}
Here, we let
\footnote{$\|\cdot\|_{M}$ is just for convenience, we don't pursue the optimal regularity required for the coefficients $\beta_{j}$ or $c$.}
\begin{equation}\label{eq: normM}
\|f\|_{M}:=\sup_{|\alpha|, |\beta|\leq 10}\|x^{\alpha}\partial^{\beta}f\|_{L_{x}^{\infty}}.
\end{equation}
It should be mentioned that,
one still has the mass conservation law
\footnote{Strictly speaking, such a mass conservation is not necessary for the construction as far as the mass is changing in a slow way, but conceptually simplifies the proof anyway.} for equation \eqref{eq: workingnls}.This is because $\|u\|_{2}=\|e^{-W}X\|_{2}=\|X\|_{2}$, and \eqref{eq: rough} admits the pathwise mass conservation law.

Recall the goal is now to construct log-log blow up solutions to \eqref{eq: workingnls}.
We will focus on the 1d case, as the 2d case follows similarly.

The observation is that, actually, a conceptually harder problem has already been studied in literature.
In \cite{raphael2006existence},\cite{raphael2009standing}, the authors studied the standing ring blow up solution for
(radial) quntic NLS on $\mathbb{R}^{d}$, $d\geq 2$,
\begin{equation}\label{eq: 6nls}
iu_{t}+\Delta u=-|u|^{4}u,
\end{equation}
and constructed (log-log type) blow up solutions which concentrate at a unit sphere $|x|=1$.
The idea is to write \eqref{eq: 6nls} in the polar coordinates,
\begin{equation}\label{eq: 6nlsmodi}
iu_{t}+\partial_{rr}u+\frac{N-1}{r}\partial_{r}u=-|u|^{4}u
\end{equation}
and view \eqref{eq: 6nlsmodi} as a mass-critical NLS on $[0,\infty)$ and
the lower order term $\frac{N-1}{r}\partial_{r}u$ as a perturbation.
Besides the bootstrap scheme explored in \cite{planchon2007existence}, the work  \cite{raphael2006existence},\cite{raphael2009standing} has two key ingredients.
\begin{itemize}
\item Since the solution is supposed to concentrate at $r\sim 1$, the term $\frac{N-1}{r}\partial_{r}u$ is supposed to behave as $\sim \partial_{r}u$, and thus subcritical.
\item Note that around origin, $\frac{\partial_{r}}{r}u$ a priori has the same  strength as $\partial_{rr}u$.
 However, the log-log blow up dynamic will concentrate its singularity around the singular point,
  and thus $u$ can be shown to be regular in a certain sense around $r=0$.
  This will make the first ingredient applicable.
\end{itemize}

From this perspective, it will  be of no surprise
that one can construct log-log blow up solutions to \eqref{eq: workingnls},
because the lower order term $\sum_{j}i\beta_{j}\partial_{j}u+cu$ is a priori subcritical,
and one does not need to handle the problem of keeping the solution regular around any specific point.

Below, we sketch the construction, highlight the key energy estimates involved in the analysis, and refer most of other  details to  \cite{raphael2006existence}, \cite{raphael2009standing}.

\subsection{Step 2: Structure of initial data and the bootstrap setup }

Let us first recall that,
there exists a continuous family of approximations  $\tilde{Q}_{b}$ of the ground state $Q$,
indexed by the parameter $b$ (see \cite{merle2003sharp}, see also \cite{merle2006sharp}, Proposition 1 and Lemma 2).
We note that $\tilde{Q}_{b}$ are uniformly smooth and localized for all $b$ small.

Following \cite{planchon2007existence}, consider the initial data
\begin{equation} \label{eq: initialdata}
u_{0}(x)=\frac{1}{\lambda_{0}^{d/2}}(\tilde{Q}_{b_{0}}+\epsilon_{0})(\frac{x}{\lambda_{0}}),
\end{equation}
where the parameters $\lambda_0$, $b_0$
and the remainder $\epsilon_0$ satisfy the orthogonality conditions\footnote{Here $\Lambda:=\frac{d}{2}+y\nabla_{y}$ is the generator of $L^{2}$ scaling.}
\begin{equation}\label{eq: modini}
\begin{aligned}
\Re (y^{2}\tilde{Q}_{b_{0}}, \overline{\epsilon_{0}})=0,\ \
\Re (y\tilde{Q}_{b_{0}}, \overline{\epsilon_{0}})=0,\ \
\Re(i\Lambda \tilde{Q}_{b_{0}}, \overline{\epsilon_{0}})=0, \ \
\Re(i\Lambda^{2}\tilde{Q}_{b_{0}}, \overline{\epsilon_{0}})=0,
\end{aligned}
\end{equation}
and the estimates
\begin{equation}\label{eq: estini}
\begin{aligned}
0<b_{0}, \quad \|\epsilon_{0}\|_{L^{2}}+b_{0}<\alpha,\\
\lambda_{0}\leq e^{-\frac{1}{\Gamma_{b_{0}}}^{4/5}},
\int |\nabla \epsilon_{0}|^{2}+|\epsilon_{0}|^{2}e^{-|y|}\leq \Gamma_{b_{0}}^{4/5}, 
\end{aligned}
\end{equation}
and
\begin{equation}\label{eq: connnn}
|E(u_{0})|\leq 1000,  |P(u_{0})|\leq 1000
\end{equation}
Here,
$\alpha$ is a small universal number,
$\Gamma_{b}\sim e^{-\pi/b}$ is a quantity frequently appearing in the log-log analysis
(see \cite{merle2004universality}),
and $E$ and $P$ denote the
energy and momentum, respectively,
i.e.,
\begin{align} \label{energy-momen}
   E(u) := \frac 12 \int |\nabla u|^2-\frac{1}{2+\frac{4}{d}}|u|^{2+\frac{4}{d}} dx, \ \
    P(u):=\Im \int \overline{u}\nabla u.
\end{align}

\begin{rem}\label{rem: furtherremark}
To make the term $i\beta_{j}\partial_{j}u+cu$ perturbative, one need to choose $\lambda_{0}$ small enough depending on the implicit constant $C$ in \eqref{eq: assumptionsforbandc}. And, by choosing $C$ large, the exceptional set we need to drop will become exponentially small with respect to $C$, though we need to choose $\lambda_{0}$ small enough to compensate large $C$.
\end{rem}

In order to construct the log-log blow up solutions to \eqref{eq: RNLS}, 
by virtue of the modulation theory \cite{weinstein1985modulational}
and the local well posedness of \eqref{eq: workingnls},
we can find some $T_{1}>0$
and parameters $b(t), \lambda(t), x(t), \gamma(t)$ such\footnote{Further standard modulation analysis plus the LWP of \eqref{eq: workingnls} show those are indeed $C^{1}$ function in $t$.} that
\begin{equation}\label{eq: gd}
u(t,x)=\frac{1}{\lambda^{d/2}(t)}(\tilde{Q}_{b}+\epsilon)(t, \frac{x-x(t)}{\lambda(t)})e^{i\gamma(t)},
\end{equation}
with the orthogonality conditions
\begin{equation}\label{eq: modulationequation}
\Re (y^{2}\tilde{Q}_{b}, \overline{\epsilon})=0, \ \
\Re (y\tilde{Q}_{b}, \overline{\epsilon})=0, \ \
\Re(i\Lambda \tilde{Q}_{b}, \overline{\epsilon})=0, \ \
\Re(i\Lambda^{2}\tilde{Q}_{b}, \overline{\epsilon})=0.
\end{equation}
and the following estimates hold:
\begin{equation}\label{eq: bootstrapassumption}
\begin{aligned}
&0<b(t), \quad \|\epsilon(t)\|_{L^{2}}+b(t)<\alpha,\\
&\forall t\leq t'\in [0,T_{1}], \quad \lambda(t')\leq \frac{3}{2}\lambda(t),\\
&\lambda(t)\leq e^{-\frac{1}{\Gamma_{b}}^{2/3}}, \\
&\int |\nabla \epsilon(t)|^{2}+|\epsilon(t)|^{2}\leq \Gamma_{b}^{2/3},\\
& t_{k+1}-t_{k}\lesssim k\lambda(t_{k})^{2}\sim k2^{-2k}.
\end{aligned}
\end{equation}
Note that $\lambda(t)$ is almost monotone in $t$, as $\forall t\leq t'\in [0,T_{1}], \quad \lambda(t')\leq \frac{3}{2}\lambda(t)$. 
Thus, one can find a partition $0=t_{0}<t_{1} <\cdots< t_{k_{0}}=T_{1}$,
such that $\lambda(t)\sim 2^{-k}$ if $t\in [t_{k}, t_{k+1}]$.

The key ingredient in the construction of log-log blow up solutions 
lies in the following bootstrap Lemma.
\begin{lem}\label{lem: boot}
Assume that $u$ solves \eqref{eq: workingnls} with the initial data $u_{0}$
and that \textcolor{black}{for some $T_{1}>0$} so that \eqref{eq: bootstrapassumption} holds.
Then, the following bootstrap estimates hold for all $t\in [0,T_{1}]$,
\begin{equation}\label{eq: bootstrapestimate}
\begin{aligned}
&0<b(t), \quad \|\epsilon(t)\|_{L^{2}}+b(t)<\alpha/2,\\
&\forall t\leq t'\in [0,T], \quad \lambda(t')\leq \frac{5}{4}\lambda(t),\\
&\lambda(t)\leq e^{-\frac{1}{\Gamma_{b}}^{3/4}},\\
&\int |\nabla \epsilon(t)|^{2}+|\epsilon(t)|^{2}e^{-|y|}\leq \Gamma_{b}^{3/4}, \\
&t_{k+1}-t_{k}\lesssim k\lambda(t_{k})^{2}\sim \sqrt{k}2^{-2k},
\end{aligned}
\end{equation}
\end{lem}
It is then left to prove the bootstrap lemma \ref{lem: boot}.

It was first observed in \cite{planchon2007existence} that the log-log dynamic admits a robust bootstrap scheme.  In our case, if we have 
energy and momentum conservation, then the proof of Lemma \ref{lem: boot} follows from\footnote{And the proof of Lemma \ref{lem: boot} is easier.} the proof of Lemma 6 in \cite{planchon2007existence}. (It will follow almost line by line up to some natural modification.)
See also the parallel bootstrap arguments in \cite{raphael2006existence}, \cite{raphael2009standing}.
Our statement of the bootstrap Lemma \ref{lem: boot} is more close to the version in \cite{colliander2009rough},
but the frequency truncated versions of energy and momentum are not required here.

One technical problem here is that, the energy and momentum are not conserved along the flow generated by equation \eqref{eq: workingnls}.

However, as observed in \cite{planchon2007existence}
and in \cite{colliander2009rough}\footnote{In \cite{colliander2009rough}, the situation is more involved because of the low regularity of the solution,
and thus the frequency truncated versions of energy and momentum need to be handled.},
under the bootstrap assumption $\lambda$ is so small in the sense that
\begin{equation}
\lambda\lesssim e^{-e^{\frac{C}{b}}}
\end{equation}
and only the size of  $\lambda^{2}E$ and $\lambda P$ will be relevant in the log-log blow up analysis.
This fact enables us to reduce to the proof of
the key bootstrap estimates in Lemma \ref{lem: boot}
to that of the  estimates of energy and momentum.

\begin{lem}\label{lem: enegyestimate}
Assume that $u$ solves \eqref{eq: workingnls} with the initial data $u_{0}$ given by \eqref{eq: initialdata}
and that estimates \eqref{eq: bootstrapassumption} hold for $t\in  [0, T_{1}]$.
Then, one has for all $t\leq T_{1}$,
\begin{equation}\label{eq: energycon}
\lambda^{2}|E(u)|\leq \Gamma_{b}^{10},
\end{equation}
and
\begin{equation}\label{eq: momentum}
\lambda |P(u)|\leq \Gamma_{b}^{10}.
\end{equation}
\end{lem}
Step $3$ below is devoted to the proof of Lemma \ref{lem: enegyestimate}.

Before we go to next step, we would like to take this chance to explain what do we mean that the term $\sum_{j}ib_{j}u+cu$ is subcritcal and what do we mean by only $\lambda^{2}E$ rather than $E$ itself is involved in the log-log  analysis.

The starting point of log-log analysis is to plug in the ansatz \eqref{eq: gd} into \eqref{eq: workingnls} and derive an equation about $\epsilon$, and apply the four orthogonality condition \eqref{eq: modulationequation} to derive the equation for $b(t), \lambda(t), x(t), \gamma(t)$. And as explained in \cite{merle2005blow}, it is more favorable to work in the re-scaled time variable $s$, such that $ds=\lambda^{-2}dt$.   And pure algebraic computation gives

\begin{equation}\label{eq: mod}
i\partial_{s}\epsilon+b_{s}\partial_{b}\tilde{Q}_{b}=L_{1}(\epsilon, b(s),x(s),\gamma(s))+L_{2}(\epsilon, b(s),x(s),\gamma(s))
\end{equation}
Here we use $L_{2}$ to denote all the extra terms caused by $\sum_{j}i\beta_{j}\partial_{j}u+cu$.
Let us also take $y=\frac{x-x(t)}{\lambda(t)}=\frac{x-x(s)}{\lambda(s)}$. Compute that
\begin{equation}\label{eq: 1}
 \frac{d}{dt}\frac{1}{\lambda^{d/2}(t)}(\tilde{Q}_{b}+\epsilon)(\frac{x-x(t)}{\lambda(t)})e^{i\gamma(t)}=\frac{1}{\lambda^{d/2+2}}(\partial_{s}\epsilon+b_{s}\partial_{b}\tilde{Q}_{b})(s,y)e^{i\gamma(t)}+\text { othe terms},
 \end{equation}
 and 
 \begin{equation}\label{eq: 2}
 \sum_{j}i\beta_{j}\partial_{j}u+cu=\frac{1}{\lambda^{d/2+1}}(t)\sum_{j}i\beta_{j}(x)\partial_{j}(\tilde{Q}_{b}+\epsilon)(y)+\frac{1}{\lambda^{d/2}(t)}c(x)(\tilde{Q}_{b}+\epsilon)(y)
 \end{equation}
 
 What one really need to do is to observe there is a $\frac{1}{\lambda^{d/2+2}}$ in \eqref{eq: 1} and there is a factor of $\frac{1}{\lambda(t)^{d/2+1}}$ in the $\beta_{j}$ related terms and a factor of $\frac{1}{\lambda^{d/2}}$ in the $c$ related terms. And a direct computation will give $L_{2}$ in \eqref{eq: mod} is 
 \begin{equation}
 L_{2}=\lambda \beta_{j}(\lambda y+x(s))\partial_{j}(\tilde{Q}_{b}+\epsilon)(y)+\lambda^{2}c(\lambda y+x(s))(\tilde{Q}_{b}+\epsilon)e^{i\gamma(t)}
 \end{equation}
 
 One crucial step in log-log analysis, \cite{merle2003sharp},\cite{merle2004universality},\cite{merle2005blow},\cite{merle2006sharp}  is to derive a local virial, which is of form 
 \begin{equation}\label{eq: lv}
 b_{s}\geq -\Gamma_{b}^{1-C\eta}-2\lambda^{2}E
 \end{equation}
 Here $C \eta \ll 1$, and we note the exact log-log law is corresponding to $b_{s}\sim -\Gamma_{b}$.

First, we can see only $\lambda^{2}E$ is invovled in the analysis, and  \eqref{eq: energycon} is enough to make this term completely perturbative.
(We remark if one further explore the modulation theory, then only $\lambda P$ will be involved, and as far as its size is neglectable compared with $\Gamma_{b}$, it will not impact the analysis.)

Second, we remark \eqref{eq: lv}, though highly nontrivial and is one of the key breakthrough  in \cite{merle2005blow},\cite{merle2003sharp}, starts from ($L^{2}$) pairing \eqref{eq: mod} with some $Q$ based well localized smooth function $\psi$, thus the extra term $L_{2}$ will only cause a perturbation of form 
\begin{equation}
O(\int \lambda|\psi||\sum_{j} i\beta_{j}(\lambda y+x(t))\partial_{j}(\tilde{Q}_{b}+\epsilon)(y)|)+O(\int \lambda^{2}|\psi||c(\lambda_{j}y+x(s))||\tilde{Q}_{b}|+|\epsilon|)(y)
\end{equation}
Thanks to the fact $\psi$ is nice and $\epsilon$ is bounded in $H^{1}$ due to \eqref{eq: bootstrapassumption}, and $\qb$ is nice, the above term is bounded by $\lambda$, which is much smaller than $\Gamma_{b}^{100}$ via \eqref{eq: bootstrapassumption}. Thus this part is also completely peturbative (or as aforementioned, subcritical).
\subsection{Step 3: Energy type estimates}\label{subsec: energy}
In this subsection, we prove Lemma \ref{lem: enegyestimate}.
Let us start with the control of energy.
We shall use the Einstein summation below.
Straightforward computations show that
\begin{equation}
\begin{aligned}
\frac{d}{dt}E(u)
=&\Re \int (- \partial_{kk}\bar{u}-|u|^{4}\bar{u})(-\beta_{j}\partial_{j}u+icu)\\
=&\Re \int  \partial_{kk}\bar{u}\beta_{j}\partial_{j}u+\Re \int |u|^{4}u\sum_{j}\partial_{j}\beta_{j}\partial_{j}u+O(\|u\|^{4}_{L_{x}^{2}}\|u\|_{H_{1}}^{2})+O(\|u\|_{H^{1}}^{2})\\
=&\Re \int  \partial_{kk}\bar{u}\beta_{j}\partial_{j}u+\Re \int |u|^{4}u\sum_{j}\partial_{j}\beta_{j}\partial_{j}u+O(\|u\|_{H^{1}}^{2})\\
=&\Re \int -\partial_{k}\bar{u}\beta_{j}\partial_{jk}u+\Re \int |u|^{4}\bar{u}\beta_{j}\partial_{j}u+O(\|u\|_{H^{1}}^{2})\\
=&\Re \int -\frac{1}{2}\partial_{j}(|\nabla u|^{2})\beta_{j}+\Re \int \frac{1}{6}\partial_{j}|u|^{6}\beta_{j}+O(\|u\|_{H^{1}}^{2}) \\
=&\Re \int \frac{1}{2}(|\nabla u|^{2})\partial_{j}\beta_{j}-\Re \int \frac{1}{6}|u|^{6}\partial_{j}\beta_{j}+O(\|u\|_{H^{1}}^{2}) \\
=&O(\|u\|_{H^{1}}^{2})
\end{aligned}
\end{equation}
(Note that there is an extra cancellation thanks to the real valued property of $\beta_{j}$.)
This is essentially formula (5.20) in \cite{SZ20}. This formula, though simple, saves one derivative  by naively plug in all the terms in. This is crucial for us to close the construction in $H^{1}$ based formula.

The point is that, in the bootstrap regime \eqref{eq: bootstrapassumption},
one has
\begin{equation}
\|u\|_{H^{1}}\lesssim \frac{1}{\lambda(t)}.
\end{equation}
Thus, we have that for any $t\in [t_{k}, t_{k+1}]\subseteq[0,T_{1}]$,
\begin{equation}
|E(u(t))-E(u_{0})|\lesssim \int_{0}^{t}\frac{1}{\lambda^{2}(\tau)}d\tau\lesssim \sum_{l\leq k+1}l\lesssim \lambda(t_{k})^{-1/10}\sim 2^{k/10}.
\end{equation}

This plus the fact $\lambda(t)\lesssim e^{-e^{C/b}}$ gives the desired result.

The control of momentum is easier, and one simply compute as
\begin{equation}
|\frac{d}{dt}P(u)|\lesssim O(\|u\|_{H^{1}}^{2}) 	
\end{equation}
To see this, just observed the $\Delta u-|u|^{4}u$ in the \eqref{eq: workingnls} preserves the momentum, and all other extra term will only cause one derivative  loss. (More careful computation can upgrade the bound to $\lesssim \|u\|_{H^{1}}$, but we don't need that there.)

The desired result follows by argue similarly as the control of energy.

\subsection{Step 4: Conclusion}

The energy estimate in Lemma \ref{lem: enegyestimate} ensures that,
though the energy and momentum are not conserved,
the proof of bootstrap Lemma in \cite{planchon2007existence} is still applicable
to our bootstrap lemma \ref{lem: boot}.
Again, we need the observation from  \cite{raphael2006existence} and \cite{raphael2009standing},
the lower order term $i\beta_{j}\partial_{j}u+c$ scales in a subcritical way.

We briefly recall the bootstrap process in \cite{planchon2007existence} for the convenience  of the readers.

The mass conservation law gives $\|u_{0}\|_{2}^{2}\sim b^{2}+\|Q\|_{2}^{2}+\|\epsilon\|_{2}^{2}$, this close the first line of \eqref{eq: bootstrapestimate} except for $b>0$.

Modulation analysis gives $-\lambda_{s}/\lambda\sim b$, thus $\lambda$ is essentially monotone.

The key is the local virial Proposition 2 in \cite{merle2003sharp} and Lyapounov control  Proposition 4 in \cite{merle2003sharp} can still be obtained in a same way  since we have the desired energy estimate Lemma \ref{lem: enegyestimate} and observation that $i\beta_{j}\partial_{j}u+c$ scales in a subcritical way, Those two estimates controls the dynamic of $b$, which decides the dynamic of $\lambda$, this gives the third line and fifth line of \eqref{eq: bootstrapestimate}.

Finally, one can see $b_{s}\sim \Gamma_{b}$, which ensures $b>0$ (before blow up), and the Lyapounove control  Proposition 4 in \cite{merle2003sharp} gives the desired monotonicity of $\int |\nabla \epsilon|^{2}+\|\epsilon\|^{2}e^{-|y|}$, which gives the fourth line of \eqref{eq: bootstrapestimate}.

Bootstrap
\footnote{Strictly speaking, to prove bootstrap Lemma \ref{lem: boot}, one needs to recover local virial and Lyponouv control which are crucial in the log-log analysis, those ingredients plus modulation theory and bootstrap Lemma \ref{lem: boot}, will imply Theorem \ref{thm: tech}.}
Lemma \ref{lem: boot} suffices to yield Theorem \ref{thm: tech}.
Most constructions of log-log blow up dynamics in literature have such an associated part and they are more or less essentially similar and are well understood now. We refer to \cite{colliander2009rough}, \cite{planchon2007existence}, \cite{raphael2006existence}, \cite{raphael2009standing}. See also \cite{fan2020construction} for a short summary.

\bibliographystyle{abbrv}
\bibliography{BG_2}
\end{document}